\documentclass{amsart}
\usepackage{times}
\usepackage{amsfonts}
\usepackage{amsmath}
\usepackage{amscd}
\usepackage{amssymb}
\usepackage{latexsym}
\usepackage{amsthm}
\usepackage[bookmarks,colorlinks,pagebackref]{hyperref} 

\theoremstyle{plain}

\swapnumbers
\newtheorem{theorem}[subsection]{Theorem}
\newtheorem{corollary}[subsection]{Corollary}
\newtheorem{lemma}[subsection]{Lemma}
\newtheorem{proposition}[subsection]{Proposition}

\theoremstyle{definition}

\newtheorem{example}[subsection]{Example}

\theoremstyle{remark}

\newtheorem{remark}[subsection]{Remark}

\swapnumbers

\newcommand\preprintnote {preprint on \myhomepage}
\newcommand\myhomepage{http://www.math.ohio-state.edu/\-\~{}schoutens}

\newcommand\NYCCT{\address{Department of Mathematics\\
NYC College of Technology\\
City University of New York\\
NY, NY 11201 (USA)}
\email{hschoutens@citytech.cuny.edu}}

\newcommand{\emptyprop}{q}

\newcommand \complet[1]{\widehat {#1}}

\newcommand \id{\mathfrak a}

\newcommand \iso{\cong}

\newcommand \map[1]{{\newcommand{\tmpprop}{#1q}  \if\tmpprop\emptyprop \to\else \xrightarrow{{\phantom{i}{#1}\phantom{i}}}\fi}} 
\newcommand \maxim{\mathfrak m}
\newcommand \nat{\mathbb N}

\newcommand \pow[2]{#1[[#2]]}

\newcommand \pr{\mathfrak p}

\newcommand \range [2]{#1,\dots,#2}
\newcommand \resp[1]{(respectively, #1)} 

\newcommand \rij[2]{(#1_1,\dots,#1_{#2})}

\newcommand \op\operatorname

\newcommand{\commdiagram}[9][]{%
\begin{equation}
{\newcommand{\tmpprop}{#1q} 
\if\tmpprop\emptyprop \relax\else \label{#1}\fi}
\begin{aligned}%
\mbox{
\begin{picture}(130,90)%
\put(120,70){\vector( 0,-1){50}}%
\put(10,80){\vector( 1, 0){100}}%
\put(0,70){\vector( 0,-1){50}}%
\put(10,10){\vector( 1, 0){100}}%
\put(115,80){\makebox(0,0)[l]{$#4$}}%
\put(5,80){\makebox(0,0)[r]{$#2$}}%
\put(115,10){\makebox(0,0)[l]{$#9$}}%
\put(5,10){\makebox(0,0)[r]{$#7$}}%
\put(-3,50){\makebox(0,0)[r]{$#5$}}
\put(123,50){\makebox(0,0)[l]{$#6$}}
\put(60,3){\makebox(0,0)[c]{$#8$}}
\put(60,88){\makebox(0,0)[c]{$#3$}}
\end{picture}}
\end{aligned}
\end{equation}}

\newcommand\commtrianglefront[7][]{%
\begin{equation}
{\newcommand{\tmpprop}{#1q} 
\if\tmpprop\emptyprop \relax\else \label{#1}\fi}
\begin{aligned}%
\mbox{
\begin{picture}(120,80)%
\put(55,68){\vector(-1,-2){30}}
\put(65,68){\vector(1,-2){30}}
\put(30,5){\vector(1,0){60}}
\put(60,75){\makebox(0,0)[c]{$#2$}}
\put(25,5){\makebox(0,0)[r]{$#4$}}
\put(95,5){\makebox(0,0)[l]{$#6$}}
\put(60,0){\makebox(0,0)[c]{$#5$}}
\put(37,43){\makebox(0,0)[r]{$#3$}}
\put(83,43){\makebox(0,0)[l]{$#7$}}
\end{picture}}
\end{aligned}
\end{equation}}

\newcommand\commtriangleback[7][]{%
\begin{equation}
{\newcommand{\tmpprop}{#1q}
\if\tmpprop\emptyprop \relax\else \label{#1}\fi}
\begin{aligned}%
\mbox{
\begin{picture}(120,80)%
\put(55,70){\vector(-1,-2){30}}
\put(65,70){\vector(1,-2){30}}
\put(30,5){\vector(1,0){60}}
\put(60,75){\makebox(0,0)[c]{$#2$}}
\put(25,5){\makebox(0,0)[r]{$#6$}}
\put(95,5){\makebox(0,0)[l]{$#4$}}
\put(60,0){\makebox(0,0)[c]{$#5$}}
\put(37,43){\makebox(0,0)[r]{$#7$}}
\put(83,43){\makebox(0,0)[l]{$#3$}}
\end{picture}}
\end{aligned}
\end{equation}}

\newcommand \ch{characteristic}
\newcommand \homo{homomorphism}
\newcommand \CM{Coh\-en-Mac\-au\-lay}
\renewcommand\iff{if, and only if,}
\newcommand \DVR{discrete valuation ring}

\newcommand \ul[1]{{#1}^*}
\newcommand \seq[2]{#1\mathstrut_{#2}}
\newcommand \BS{Brian\c{c}on-Skoda}
\newcommand \usc[1]{\op{usc}(#1)}
\newcommand \ctc[1]{\op{ctc}(#1)}
\newcommand \tc[1]{\op{tc}(#1)}
\renewcommand \sc[1]{\op{sc}(#1)}
\newcommand \los{\L os' Theorem}
\renewcommand \inf[1]{\op{Inf}(#1)}
\newcommand \sep[1]{#1_{\text{sep}}}

\title {The asymptotic behavior of solid closure in mixed characteristic}
\author{Hans Schoutens}
\date{\today}
\NYCCT
\thanks{Partially supported by a PSC-CUNY grant}

\begin{document}

\begin{abstract}
We study how solid closure in mixed \ch\ behaves  after taking ultraproducts. The   ultraproduct 
will be chosen so that we land in equal \ch, and therefore can make a
comparison with tight closure.
 As a corollary we get an asymptotic version
of the Hochster-Roberts invariant theorem in dimension three: if $R$ is a
mixed \ch\ (cyclically) pure $3$-dimensional local subring of a regular local ring $S$, then
$R$ is \CM, provided the ramification of $S$ is large with respect to its
dimension and residual \ch, and with respect to the multiplicity of $R$.
\end{abstract}

\maketitle

\section{Ultra- versus cata-}

Solid closure was introduced by Hochster in \cite{HoSol,HoTCSol} as a
potential substitute for  tight closure  in mixed \ch.  In this note, we
comment on some of its properties, but as the title indicates, only  
`asymptotically', that is to say, after taking an ultraproduct (see
\S\ref{s:asym} for an elaboration on the terminology).
More precisely, let $\seq Aw$ be a sequence of (commutative) rings (with
identity), indexed by an infinite index set endowed
with a non-principal ultrafilter, which, for technical reasons, we also assume
to be countably incomplete.\footnote{Suffice it here to say that this   set-theoretic notion 
can be realized on any infinite index set and  holds automatically when the index
set is countable. Moreover, it is consistent with ZFC (='usual set theory') that
every ultrafilter is
countably incomplete.} The ultraproduct of the $\seq Aw$ is again a ring
$A$, realized as the quotient of the product of the $\seq Aw$ modulo the ideal
of all sequences almost all of whose entries are zero (with \emph{almost all}
one means for all indices in some member of the ultrafilter).
We sometimes refer
to the $\seq Aw$ as \emph{components} of $A$, although they are not uniquely
defined by $A$. If $\mathcal P$ is a property of rings, then we say that $A$
has property \emph{ultra-$\mathcal P$} if almost all $\seq Aw$ have
property $\mathcal P$. If a property $\mathcal P$ is first-order, then
ultra-$\mathcal P$ is the same as $\mathcal P$ by \los. For instance, being
local is a first-order property,\footnote{We will call a ring $R$ \emph{local}
if it has a unique maximal ideal $\maxim$, and we denote this by $(R,\maxim)$
(in the literature one sometimes uses the term \emph{quasi-local} in the non-Noetherian
case). Note that $R$ is local \iff\ the sum of
any two non-units is a non-unit, indeed a first-order property.} so that ultra-local
is the same as local. However, most properties are not
first-order (mostly because they require quantification over ideals or involve
infinitely many statements). For instance, an \emph{ultra-Noetherian} local ring is an
ultraproduct of Noetherian local rings, and in general is 
no longer Noetherian (in fact, its prime spectrum is infinite and can be quite
complicated;   for some instances of this, see \cite{OlbSay,OlbSaySha,OlbSha}).
We
will only be concerned with a certain subclass of ultra-Noetherian local rings,
those of finite embedding dimension. An ultra-Noetherian local ring has
embedding dimension $n$ \iff\ almost all of its components have embedding
dimension $n$ (because having embedding dimension $n$ is a first-order
property). Ultra-Noetherian local rings of finite
embedding dimension already appeared as an essential tool in the earlier papers on non-standard tight closure
(\cite{SchAsc,SchLogTerm,SchNSTC,SchBCM}), and were used in
\cite{SchMixBCM,SchMixBCMCR} to get some asymptotic versions of the homological
conjectures in mixed \ch. In \cite{SchFinEmb} a systematic study of this class
will be carried out, leading to some improved asymptotic versions of the
intersection theorems in mixed \ch. The present note does not require the full
development of this theory, and we will review whatever we need.  

In the latter papers, an important technique to study local rings
of finite
embedding dimension is through their completion, since this is always Noetherian.
This leads to a second variant of a property $\mathcal P$: we call a local ring $R$ of finite
embedding dimension  \emph{cata-$\mathcal P$} if its completion has property
$\mathcal P$.\footnote{Because of its versatility, I opted to introduce the
prefix \emph{cata-} instead of using the more traditional adverbially
constructions   \emph{analytically} of \emph{formally}; similarly,
the prefix \emph{ultra-} replaces terms such as \emph{non-standard} or
\emph{generic} from the older papers.}
In case of an ultra-Noetherian local ring $(R,\maxim)$ of finite embedding
dimension, because of    saturation properties of ultraproducts, its completion equals its
separated quotient $\sep R$ (see for instance \cite[Lemma 5.1]{SchFinEmb}) defined
as the homomorphic image of $R$ modulo its
\emph{ideal of infinitesimals} $\inf R:=\bigcap_n\maxim ^n$; we call $\sep R$
the \emph{separated ultraproduct} of the components $\seq Rw$. (The   term
`cata' was chosen because of this fact.) The following example of the close connection between
 the ultra-variant of a property and its cata-variant was already observed in \cite[Corollary
1.14]{Yo}.

\begin{theorem}\label{T:ureg}
An ultra-regular local ring of finite embedding dimension is cata-regular.
\end{theorem}
\begin{proof}
Suppose $(R,\maxim)$ is the ultraproduct of regular local rings $(\seq
Rw,\seq\maxim w)$ of dimension $d$ and let $\seq{\mathbf
x}w:=(x_{1w},\dots,x_{dw})$
be a regular system of parameters in $\seq Rw$. The \emph{ultraproduct} $x_i$ of
the
$x_{iw}$ (that is to say, the image of the sequence $(x_{iw}|w)$ in $R$) gives
rise to a $d$-tuple $\mathbf x:=\rij xd$
 in $R$ generating $\maxim$, whence generating $\maxim\sep R$.
So remains to show that $\sep R$ has dimension $d$. By Krull's principal ideal
theory, its dimension is at most $d$, so suppose towards a contradiction that
it were less. Hence after renumbering, $x_d^n\in\rij x{d-1}\sep R$, for some
$n$. Contracting back to $R$, we get that $x_d^n$ lies in $\rij x{d-1}R+\inf
R\subseteq (x_1,\dots,x_{d-1},x_d^{n+1})R$. Writing out   the latter relation,
we get 
that  $x_d^n\in\rij x{d-1}R$. By
\los, $x_{dw}^n\in (x_{1w},\dots,x_{d-1,w})\seq Rw$ for almost all $w$,
contradiction.
\end{proof}

From the proof it also follows that $R$ has the same \emph{ultra-dimension}
(=dimension of almost all of its components) as \emph{cata-dimension}
(=dimension of its completion).
With this additional assumption,
the converse of Theorem~\ref{T:ureg} also holds. This is explained in
\cite[Theorem 8.1]{SchFinEmb}, where it is also shown that an ultra-Noetherian
local ring 
has the same ultra-dimension as cata-dimension 
 \iff\
the parameter degree of its components is bounded. Recall
that the \emph{parameter degree} of a Noetherian local ring $A$  is the
least co-length of a parameter ideal of $A$ (the \emph{co-length} of an ideal
$\id$ is the length of $A/\id$; a \emph{parameter ideal} is an ideal generated
by a system of parameters). Multiplicity is a lower bound for parameter degree
(with equality \iff\ the ring is \CM, provided the residue field is infinite; see \cite[Corollary 3.3]{SchABCM}).


\section{Closure operations}\label{s:ch}

In this section $R$ is an ultra-Noetherian local ring of finite embedding
dimension, say realized as the ultraproduct of Noetherian local rings $\seq
Rw$ of bounded embedding dimension.  We want to
introduce two closure operations on $R$: ultra-solid closure and cata-tight
closure.

\subsection{Ultra-solid closure}
For the definition of solid closure, see \cite{HoSol}. To maintain a consistent
notational scheme, we will write $\sc \id$ for the solid closure of an ideal
$\id$, rather than Hochster's $\id^\star$.   We say that $z\in R$  lies in
the \emph{ultra-solid closure} $\usc I$  of an ideal $I\subseteq R$, if we can find $\seq
zw\in \seq Rw$ and $\seq \id w\subseteq \seq Rw$ with respective ultraproducts
$z$ and $\id$, such that $\id\subseteq I$ and $\seq zw\in\sc{\seq \id w}$, for
almost all $w$.   In case $I$ is an \emph{ultra-ideal} (sometimes called an
\emph{induced} ideal), that is to say, is itself an ultraproduct
of ideals $\seq Iw$, then $z$ lies in its ultra-solid closure \iff\
almost all $\seq zw$ lie in the solid closure of $\seq Iw$. In other words, the
ultra-solid closure of an ultra-ideal is the ultraproduct of the solid
closures  of its components (and so again an ultra-ideal). Note that a finitely
generated ideal, whence in particular an $\maxim$-primary ideal, is an ultra-ideal.

\subsection{Cata-tight closure}
To define cata-tight closure, which will be derived from tight closure, we need to
make an assumption on the \ch, namely that $R$ is \emph{cata-equi\ch}, meaning that
its completion (or equivalently, its separated quotient) has equal \ch.   
Although in positive
\ch\
the more common notation for the tight closure 
of an ideal   is $\id^*$,   we will use instead, in either   \ch, the
notation $\tc\id$. For the next definition, we assume that $R$ is
cata-equi\ch, but there is no not need for assuming it is ultra-Noetherian--having finite embedding dimension suffices. We say that $z$ lies in
the
\emph{cata-tight closure} $\ctc I$ of $I$, if the image of $z$ in $\complet R$
lies in
the tight closure of $I\complet R$. In other words, $\ctc I=\tc{I\complet R}\cap R$.
We have a choice in picking a tight closure
operation in equal \ch\ zero: there are the  'classical' ones introduced by
Hochster and Huneke in \cite{HHZero}, and there are the `non-standard' variants
defined in \cite{SchAsc,SchNSTC}. As a rule, we will use \emph{generic tight closure} as
defined in the latter papers (see also the next paragraph). It has all the
properties we want it to have:
it is trivial on regular local rings, it `captures colons' and it is `persistent'
(for details see \cite[\S6.19]{SchAsc}). 

\subsection{More closure operations}\label{s:cl}
Of course nothing prevents us
from
switching around these definitions and introduce also `cata-solid closure' and
`ultra-tight closure'. We only would be able to say  something sensible about the
former when it actually coincides with cata-tight closure (namely, when $\complet R$ has
positive \ch); as for the latter, it is very akin to generic tight closure in
the main case where the components have positive \ch\ but $R$ has zero \ch,
and so will add nothing new. There is one more closure operation in equal \ch\
which is even smaller than tight closure (but conjecturally coincides with
it), to wit, \emph{plus closure}. In \ch\ $p$, it is defined for a Noetherian 
local domain $A$ as contraction from the \emph{absolute integral closure} $A^+$ (=integral
closure of $A$ in an algebraic closure of its field of fractions; recall
that if $A$ is moreover excellent, then $A^+$ is a big \CM\ $A$-algebra by
\cite{HHbigCM}).
In equal \ch\ zero, it is defined by contraction from a canonically defined big
\CM\
$A$-algebra
  $\mathfrak B(A)$ (see \cite{SchBCM} for the affine case and \cite{SchAsc}
for the general). On occasion, we will thus encounter \emph{cata-plus closure} as
the pre-image of plus closure in $\complet R$.

We say that an ideal $I\subseteq R$ is \emph{ultra-solidly closed}
\resp{\emph{cata-tightly closed}}
if $I$ is equal to its ultra-solid closure \resp{cata-tight closure}. We cannot expect
cata-tight closure to be a true `tight' closure operation, since it always contains
the $\maxim$-adic closure of the ideal (see Lemma~\ref{L:clos} for a
description of the closure of an ideal). However, $\maxim$-primary ideals are
$\maxim$-adically 
closed and hence their cata-tight closure will be the most accessible.

\subsection{From mixed to equal \ch}
We now turn to the issue of enforcing equal \ch\ for the separated
ultraproduct, starting from mixed \ch. One way is by letting the components
$\seq Rw$ have unbounded residual \ch, that is to say, for each $n$, the \ch\
of the residue field of $\seq Rw$ is $\geq n$ for almost all $w$. This has
a consequence that the ultraproduct $R$ has residual \ch\ zero, whence so does
$\sep R=\complet R$. However, there is a second way to get an equal \ch\
separated ultraproduct:

Let $(A,\pr)$ be a local ring of  residual  \ch\ $p$.
We define the \emph{ramification index} of $A$,   as the
supremum of all $n$ for which $p\in\pr^n$. We call $A$ \emph{unramified}, 
 if its ramification index is one, and   \emph{infinitely
ramified}, if it is infinite, that is to say, if $p\in\inf A$. (Note: if the
residual \ch\ is zero, then we will call $A$ also \emph{unramified}). 
 If $A$ is Noetherian, or
more generally, separated,
and infinitely ramified, then in fact it has  equal \ch\ $p$ (in the literature
this situation is erroneously referred to as `unramified'). However, in
  general, a local ring can have \ch\ zero and be infinitely ramified:
if $\seq Rw$ are mixed \ch\ Noetherian local rings of residual \ch\ $p$ and
unbounded ramification index (in the sense that for each $n$, almost all $\seq
Rw$ have ramification index  $\geq n$), then their
ultraproduct $R$
has
\ch\
zero
and is infinitely ramified. In particular, the separated ultraproduct $\sep R$
of the $\seq Rw$ has equal \ch\ $p$.  

In summary, $R$ is cata-equi\ch, if   either $R$
itself has equal \ch, or otherwise, is infinitely ramified. Unfortunately,
solid closure does not behave that well in equal \ch\ zero (see
Example~\ref{E:sol0}),
so
that we can only compare our two new closure operators when the completion has
prime \ch.

\begin{proposition}\label{P:uscc}
If $R$ is an ultra-Noetherian local domain of finite
embedding dimension, which has prime \ch\   or is infinitely ramified, then
 $\usc I\subseteq\ctc I$, for all   ideals $I\subseteq R$.
\end{proposition}
\begin{proof}
Let $z\in\usc I$ and let $\seq zw\in\seq Rw$ and $\seq \id w\subseteq \seq Rw$
have
respective ultraproducts $z$ and $\id\subseteq I$ with almost all $\seq zw$ in
the solid closure of $\seq \id w$. By
\cite[Proposition 5.3]{HoSol},
there exists for almost all $w$ a solid $\seq Rw$-algebra $\seq Sw$ such that
$\seq
zw\in\seq \id w\seq Sw$. Recall that $\seq Sw$ is \emph{solid} when there
exists a non-zero $\seq Rw$-linear
map $\seq\phi w\colon \seq Sw\to \seq Rw$. Moreover,  we may choose $\seq\phi
w$ so that $\seq\phi w(1)\neq 0$.  
Let $S$ be the ultraproduct of the $\seq Sw$. The ultraproduct
$\phi$ of the  $\seq\phi w$ is an $R$-linear map $S\to R$, showing that $S$
is solid as an $R$-algebra.
Let $\tilde S:=S/\inf RS$. Hence $\phi$ induces by base change an 
$\sep R$-linear map $\tilde S\to \sep R$ showing that $\tilde S$ is a solid $\sep
R$-algebra. By assumption, $\sep R$ has \ch\ $p>0$. By \los, $z\in \id\tilde
S$. Applying Frobenius to this equation   and then applying $\phi$, we get  
$cz^q\in \id^{[q]}\sep R$, for all powers $q$ of $p$, with $c:=\phi(1)\neq 0$.
Hence $z\in\tc {\id\sep R}$ and therefore  $z\in \ctc \id\subseteq \ctc I$
(recall that $\sep R=\complet R$). 
\end{proof}

\begin{example}\label{E:sol0}
The above inclusion does not hold in general in equal \ch\ zero: in
\cite{RobSol} Roberts shows that $f:=X^2Y^2Z^2$ lies in the solid closure of
$I:=(X^3,Y^3,Z^3)A$ where $A:=\pow K{X,Y,Z}$ with $K$ a field of \ch\ zero. Let
$R$ be the ultrapower of $A$ (an \emph{ultrapower} is an ultraproduct in which
all components are the same).
It follows that $f\in\usc {IR}$. On the other
hand, $\sep R\iso\pow{\ul K}{X,Y,Z}$ where $\ul K$ is the ultrapower of $K$,
so that $I\sep R$ is tightly closed and hence $IR$ is cata-tightly closed. I do
not know
of any example of an equal \ch\ zero ultra-Noetherian local domain whose components
have mixed \ch, but for which the above inclusion does not hold. There is
also the hope that the above result holds without any restriction on the \ch\ when
we replace solid closure by
\emph{parasolid closure} (promising in that respect is  \cite[Theorem
4.1]{BreSol} showing that every ideal in a regular local ring is parasolidly
closed).
\end{example}

It is an interesting question whether we can have equality in
Proposition~\ref{P:uscc}. However,
since an ultra-ideal cannot be $\maxim$-adically closed unless it is $\maxim$-primary, whereas
cata-tight closures are always $\maxim$-adically closed, we can only expect
equality for $\maxim$-primary ideals. Using a result of Smith, we can prove this in a
special case. Let $(R,\maxim)$ be a local ring of finite embedding dimension. A  tuple $\mathbf x$
in $R$ is called a \emph{system of cata-parameters} (or a \emph{generic
tuple}), if its image in $\sep R$ is a system of parameters. Any ideal
generated by a system of cata-parameters will be called a
\emph{cata-parameter} ideal. In particular, a cata-parameter ideal $I$
is $\maxim$-primary and $I\sep R$ is a parameter ideal. 


 \begin{theorem}\label{T:par}
Let $R$ be an   ultra-Noetherian local ring  of finite embedding dimension, 
whose separated quotient is an equidimensional, prime \ch\ reduced local ring. If
$I$ is a cata-parameter ideal in $R$, then $\usc I=\ctc I$.
\end{theorem}
\begin{proof}
One inclusion is clear from Proposition~\ref{P:uscc}. Hence assume that $z$
lies in $\ctc I$, whence its image in $\tilde R:=\sep R$ lies in $\tc{I\tilde R}$.
Let $\pr$ be  a minimal prime  of $\tilde R$. Since $\tilde R$ is  equidimensional,
$I(\tilde R/\pr)$ is a parameter
ideal. By persistence, $z$ lies in $\tc{I(\tilde R/\pr)}$, whence in $I(\tilde R/\pr)^+$
by \cite{SmParId},  where
$(\tilde R/\pr)^+$ is the  absolute integral closure  of the complete local
domain $\tilde R/\pr$ (see \S\ref{s:cl}). It follows that there
exists a finite extension $\tilde R/\pr\subseteq \tilde S(\pr)$ of local domains
such that
$z\in I\tilde S(\pr)$. Let $\tilde S$ be the direct sum of all $\tilde S(\pr)$,
where $\pr$ runs over all minimal primes of $\tilde R$. Since $\tilde
R$ is reduced, the natural map $\tilde R\to \tilde S$ is finite and injective.
Hence we can lift this to a finite extension $R\subseteq S$, such that $\tilde
S\iso S/\inf RS$. From $z\in I\tilde S$ and the fact that $\inf R\subseteq I$,
we get $z\in IS$. 

Choose finite local extensions $\seq Rw\subseteq
\seq Sw$ whose ultraproduct is $R\subseteq S$. Let $\seq zw$ and $\seq Iw$ be
such that their ultraproducts are $z$ and $I$ respectively. By \los, almost
each $\seq zw$ lies in $\seq Iw\seq Sw$ whence in the solid
closure of $\seq Iw$ since finite extensions are formally solid by \cite[Remark
1.3]{HoSol}. In conclusion, we showed that $z\in \usc I$.
\end{proof}

%
%
%

The argument in the proof actually shows that the cata-plus closure (see
\S\ref{s:cl}) is always contained in the ultra-solid closure for
$\maxim$-primary ideals, for $R$ as in the statement. Therefore, if plus-closure
equals tight closure
in positive \ch, ultra-solid closure, cata-tight closure and cata-plus closure
are all the same
on $\maxim$-primary ideals. If $\sep R$ has equal \ch\ zero, the above argument
does not work since $\mathfrak B(\sep R)$ is no longer integral over $\sep R$
and hence it is not clear how to `descend' to the components.

\section{Properties of cata-tight closure}

In this section, $(R,\maxim)$ denotes a cata-equi\ch\ local ring
of finite embedding dimension. Our goal is to derive some elementary properties
of cata-tight closure on $R$. We already observed that the $\maxim$-adic closure
is contained in the cata-tight closure.

\begin{lemma}\label{L:clos}
The $\maxim$-adic closure of an ideal $I\subseteq R$ is equal to $I+\inf R$.
\end{lemma}
\begin{proof}
Since $\sep R$ is Noetherian, the $\maxim$-adic closure of $I\sep R$ is equal to
$I$ by
Krull's Intersection theorem, and the assertion follows. (No assumption on the
\ch\ is needed for this lemma).
\end{proof}

 Since in a regular local ring, every  ideal is tightly closed, we get
immediately:

\begin{proposition}\label{P:ccreg}
If $R$ is cata-regular and $I$ an ideal in $R$, then the cata-tight closure of $I$ is
equal to its $\maxim$-adic closure, that is to say, $\ctc I=I+\inf R$. In particular, any
$\maxim$-primary ideal is cata-closed.
\end{proposition}

Colon Capturing and persistence of (generic) tight closure in equal \ch\ leads
immediately to the analogous results for cata-tight closure: 
 
 \begin{proposition}\label{P:CC}
If $\rij xd$ is a system of cata-parameters in $R$, then for each $i\leq d$, we
have an inclusion $(\rij x{i-1}R:x_i)\subseteq \ctc{\rij x{i-1}R}$ (`Colon
Capturing').

If $R\to S$ is a local \homo\ and $I\subseteq R$ an ideal, then $\ctc
IS\subseteq \ctc {IS}$ (`Persistence').
\end{proposition}

\begin{remark}\label{R:CC}
There are actually many stronger versions of Colon Capturing, of which I only will
mention   one: given integers
$0\leq a_i<b_i$ for $i=\range 1d$, and a system of cata-parameters $\rij
xd$, we have an inclusion
\begin{equation}\label{eq:sCC}
(\ctc{(x_1^{b_1},\dots,x_d^{b_d})R} : x_1^{a_1}\cdots x_d^{a_d})\subseteq
\ctc{(x_1^{b_1-a_1},\dots,x_d^{b_d-a_d})R}.
\end{equation}
In positive \ch, this inclusion follows from the tight closure version 
of \eqref{eq:sCC} proven in \cite[Theorem 9.2]{HuTC}. The latter
then also gives the corresponding result in zero \ch\ by the techniques of
\cite{SchAsc}.
\end{remark}

In the next
result, we have written $\bar I$ to denote the \emph{integral closure} of an
ideal $I$.

\begin{theorem}[\BS]\label{T:cBS}
In a cata-equi\ch\  local ring $(R,\maxim)$ of finite embedding dimension, we
have for each ideal
$I\subseteq R$ an inclusion $\overline {I^d}\subseteq \ctc I$, where $d$ is the
cata-dimension of $R$.

In particular, if $R$ is cata-regular and $I$ is $\maxim$-primary, then
$\overline {I^d}\subseteq I$.
\end{theorem}
\begin{proof}
Let $\id:=I\complet R$. Clearly, $\bar I\complet R\subseteq \bar\id$. By the tight
closure \BS\ theorem, $\overline {\id^d}\subset \tc\id$, so that the
first assertion is clear. The second assertion then follows from
Proposition~\ref{P:ccreg}.
\end{proof}

 \begin{remark}\label{R:BS}
In fact, the usual \BS\ theorem gives the following stronger result: if $h$
denotes the minimal number of generators of $I\sep
R$, then we have for all $k$ an inclusion $\overline {I^{k+h}}\subseteq \ctc 
{I^{k+1}}$. Using an improvement by Aberbach and Huneke in equal \ch\ in
\cite{AbHuBS}, we   actually get the following. Assume $R$ is
cata-regular, with infinite residue field,
$I$ is $\maxim$-primary and $J\subseteq I$ is a minimal reduction of $I$, then
\begin{equation}\label{eq:AH}
\overline {I^{k+d}}\subseteq   J^{k+1}\id,
\end{equation}
for all $k$, where $\id$ is maximal among all ideals for which $\id J=\id I$
(note that $\id\subseteq (J:I)$ and hence is a proper ideal unless $I$ is its
own minimal reduction).
\end{remark}
 
\subsection{Asymptotic properties}\label{s:asym}
The next type of result explains better the
term `asymptotic' from the title:
a property (often of homological nature) holds `asymptotically' when the \ch\ (or
the ramification) is large with respect to the other
data (in a sense that has to be made more precise).  In
\cite{SchMixBCM,SchMixBCMCR}, the lower bound for the \ch\ depended on
the degrees of the polynomials defining the data. In \cite{SchFinEmb} an
improved bound for the intersection theorems will be given only  depending  on
some (more natural) invariants of the
data (like dimension and parameter degree). Proofs of these types of results
all follow the same outline: if there are counterexamples, their ultraproduct
  violates  the corresponding ultra-version, which holds because its cata-version holds 
since we are now in equal \ch.  The first assertion in the next result is just
included as an example, for it follows already from the general \BS\ theorem of
Lipman and Sathaye in \cite{LS}, which holds for all regular local rings,
regardless of their \ch. 

\begin{theorem}[Asymptotic \BS\ in mixed \ch]\label{T:aBS}
For each pair $(d,l)$, there exists a bound $B:=B(d,l)$ with the property that
if $S$ is a $d$-dimensional mixed \ch\ regular local ring   and if $I\subseteq
S$ is an ideal of co-length at most $l$, then $\overline {I^d}\subseteq I$,
 provided the residual \ch\ of $S$, or its ramification index,
  is at least $B$.
  
 In fact, under these assumptions, we have an inclusion $\overline {I^d}\subseteq
\id I$, provided $I$ has a reduction $J$ of co-length at most $l$ and 
$\id$ satisfies $\id J=\id I$.
\end{theorem}
\begin{proof}
I will only give the details   for the first assertion, following
the outline just mentioned. As for the second assertion, it follows along the
same lines, using \eqref{eq:AH} from Remark~\ref{R:BS} instead. Suppose the first
assertion is false for some pair $(d,l)$. This means
that for each $w\in\nat$, we can find  a mixed \ch\ $d$-dimensional regular local
ring $\seq Rw$  whose residual \ch\ or ramification index   $\geq w$ and an ideal
$\seq I w\subseteq
\seq Rw$
 of co-length at
most $l$, such that  $\overline {\seq I w^d}$ is not contained in $\seq  I
w$. Let $I$ and $R$ be the respective ultraproducts of $\seq I w$ and
$\seq Rw$. It follows that $R$ is cata-equi\ch\ and ultra-regular, whence
cata-regular,
by Theorem~\ref{T:ureg}. By \los, $I$ has co-length at most $l$ and
 does not contain $\overline{I^d}$, contradicting Theorem~\ref{T:cBS} (the
`cata-\BS\ theorem'). 
\end{proof}

\begin{remark}
By the same argument and under the same assumptions, there is also a bound
$B':=B'(d,l,m)$ such that
\eqref{eq:AH} holds for all  $k\leq m$ whenever the residual \ch\ or the
ramification index $\geq B'$.
\end{remark}

\section{Properties of ultra-solid closure}
Although we know very little about solid closure in mixed \ch, it is clear
that the inclusion from Proposition~\ref{P:uscc} together with the results
from the previous section should  tell us a whole lot more about its ultra-variant.
Consequently, we may hope to infer some `asymptotic' properties of solid
closure itself, at least for $\maxim$-primary ideals of bounded co-length.

Combining Propositions~\ref{P:uscc} and \ref{P:ccreg} with
Theorem~\ref{T:ureg} yields immediately:

\begin{corollary}\label{C:screg}
If $(R,\maxim)$ is an ultra-regular local ring of finite
embedding dimension which has prime \ch\ or is  infinitely ramified, then every $\maxim$-primary
ideal is ultra-solidly closed.
\end{corollary}

\begin{corollary}\label{C:usHR}
Let be $(R,\maxim)\to (S,\mathfrak n)$ be a local \homo\ of ultra-Noetherian
local rings of finite embedding dimension. If  $S$ is ultra-regular and 
has prime \ch\ or is infinitely ramified, then
$\usc \id\subseteq \id S\cap R$, for every $\maxim$-primary ideal
$\id\subseteq R$.
\end{corollary}
\begin{proof}
Let  $z\in \usc \id$. Since solid closure is persistent,
so is ultra-solid closure, and hence $z$ lies in $\usc {\id S}$, whence in
$\id S$ by Corollary~\ref{C:screg}.
\end{proof}

 Following Hochster, we call   a Noetherian local
ring  \emph{weakly S-regular} \resp{\emph{S-rational}} if every ideal
\resp{every parameter ideal} is solidly closed. For tight closure, if a single parameter ideal   is tightly
closed, then so is any other parameter ideal, but not so for solid closure.
We therefore will call $R$ \emph{weakly S-rational} if it admits at least one
solidly closed parameter ideal. The
next result shows that this notion serves some purpose.

\begin{proposition}\label{P:wSratCM}
If an analytically irreducible Noetherian local ring is weakly S-rational and
admits a big \CM\ algebra, then it is \CM.
\end{proposition}
\begin{proof}
Let $R$ be a weakly S-rational Noetherian local ring and let $B$ be a
big \CM\ $R$-algebra. By \cite[Corollary 8.5.3]{BH}, we may replace $B$ by its
$\maxim$-adic completion and assume that it is even a \emph{balanced} big \CM\
algebra.  By the argument in \cite[Proposition 7.9(c)]{HoSol} we may then
replace $R$ by its completion.  Let $\mathbf x:=\rij xd$ be a system of
parameters generating a solidly closed ideal.  Let
$I_i:=\rij xiR$. Before we prove the proposition, let us show by downward
induction on $i\leq d$ that $I_i=I_iB\cap R$. The case $i=d$ follows from our
assumption that $I_d$ is solidly closed and the fact that a big \CM\ algebra
over a Noetherian local domain is solid. Hence let $i<d$ and assume
$I_{i+1}=I_{i+1}B\cap R$.  Let $z$ be an element of $J_i:=I_iB\cap R$. By our induction
hypothesis, $z\in I_{i+1}$, say $z=y+ax_{i+1}$ with $y\in I_i$ and $a\in R$. From
$ax_{i+1}=z-y\in I_iB$ and the fact that  $\mathbf x$ is
$B$-regular,   we get $a$ lies in $I_iB$ whence in $J_i$. In conclusion, we
showed that $J_i=I_i+x_{i+1}J_i$, so that by Nakayama's lemma,   $J_i=I_i$, as
claimed. 

To complete the proof, we must show that $\mathbf x$ is $R$-regular.
To this end, suppose $zx_{i+1}\in I_i$. Since $\mathbf x$ is $B$-regular, $z$ 
lies in $I_iB$ whence in $I_i$, by our previous remark.
\end{proof}

By a standard argument (see for instance \cite[Proposition 5.6]{SchBCM}),
every ideal generated by part of a system of parameters is then contracted from
$B$. However, this does not yet prove that $R$ is S-rational. 

Again we can derive some asymptotic versions in
mixed \ch\ of the previous results. However, in view of the
restriction on the \ch\ imposed by Proposition~\ref{P:uscc} (namely that its
separated quotient have equal \ch\ $p$), we only get an asymptotic version for
large ramification index. In our first application, even the case $R=S$  leads
to new results (although trivially weakly
S-rational,  a regular local ring of
mixed \ch\ is only conjecturally weakly S-regular or even  S-rational):

\begin{theorem}\label{T:wSrat}
For each triple $(p,n,l)$ with $p$ a prime number, there exists a bound
$N:=N(p,n,l)$ with the following property.  Let $R\to S$ be a cyclically pure   \homo\ of
Noetherian  local rings of residual \ch\ $p$ and embedding dimension at most
$n$, with $S$ regular. Let $\id$ be an ideal in $R$ of co-length at most $l$.

If $S$ has ramification index at least $N$, then $\id$ is solidly closed. In
particular, if $R$ has 
moreover parameter degree at most $l$, then it is weakly S-rational.
\end{theorem}
\begin{proof}
We only need to show the first assertion, so suppose it  is false for some
triple $(p,n,l)$. Hence, there exists for each $w$, a cyclically pure   \homo\ of
Noetherian 
local rings $\seq Rw\to \seq Sw$ of residual \ch\ $p$ and embedding dimension
at most $n$,
with $\seq Sw$
 regular
of ramification index  at least $w$, and an ideal $\seq\id w\subseteq\seq Rw$ of co-length at most $l$
which is not solidly closed. Let $\id$, $R$ and
$S$ be the respective ultraproducts of the $\seq \id w$, $\seq Rw$ and
$\seq Sw$. Therefore, $S$ is infinitely ramified  and  cata-regular. Since
$\seq\id w=\seq\id w\seq Sw\cap\seq Rw$, \los\ yields
that $\id=\id S\cap R$. Moreover, $\id$ has co-length at most $l$ whence in 
particular is $\maxim$-primary. Therefore $\id$ is ultra-solidly closed by
Corollary~\ref{C:usHR}, and hence   almost all $\seq \id w$ are solidly
closed, contradiction.
\end{proof}

\begin{remark}
In his lists of open problems (\cite[Question 20]{HoTCSol}), Hochster asks--in
the hope of getting a negative answer--the following:
does $\pi^2X^2Y^2$ belong to the solid closure of the ideal $(\pi^3,X^3,Y^3)R$,
where $R:=\pow V{X,Y}$ and $V$ is a mixed \ch\ \DVR\ with  
uniformizing parameter $\pi$ and valuation $v$? According to our previous
result, for   fixed residual \ch\ $p$,  
  the answer is indeed no, provided $v(p)$ is
sufficiently large. Ironically, Hochster
asked this for $V$ an unramified \DVR\ (in fact, for $V$ equal to the ring of
$p$-adic integers), expecting this to be the easiest case to settle, yet, this
case remains open.
\end{remark}

\begin{theorem}[Asymptotic Hochster-Roberts in dimension $3$]
For each triple $(p,d,l)$ with $p$ a prime number, there exists a bound
$\rho:=\rho(p,d,l)$
with the
following property. Let $R\to S$ be a cyclically pure   \homo\ of Noetherian 
local rings of residual \ch\ $p$. Assume $S$ is regular and has dimension at
most $d$, and assume $R$ has dimension
at most three and parameter degree at most $l$. If the ramification index of
$S$ is at least $\rho$, then $R$ is \CM.
\end{theorem}
\begin{proof}
 By Theorem~\ref{T:wSrat}, if the ramification index of $S$ is at least
$N(p,d,l)$,  then $R$ is weakly
S-rational.  Since $R$ admits a big \CM\ algebra by
\cite{HoBCM3}, the result follows from Proposition~\ref{P:wSratCM} (note that
a cyclically pure subring of a regular local ring is analytically irreducible).
\end{proof}

One expects that this result is true without any restriction on the
ramification index of $S$. To derive this more general result directly from
Hochster's result on the existence of big \CM\ algebras in dimension three,
one would need to show that the big \CM\ $B$ admits an $R$-algebra \homo\ into
a big
\CM\
algebra for
$S$ (equivalently, into some faithfully flat extension of $S$)  and this is only
known if also $S$ has dimension at most three. In higher
dimensions, the existence of big \CM\ algebras is still open, and we have to  
settle for the following much weaker 
result. Recall that a tuple $\rij xn$ is called \emph{independent} (\emph{in the sense of
Lech}), if any relation $a_1x_1+\dots+a_dx_d=0$ in $R$ implies that all $a_i$
lie in the ideal $I$ generated by $\rij xn$ (equivalently, if $I/I^2$ is a
free $R/I$-module of rank $n$). Any regular sequence is independent, and
conversely, if $(x_1^t,\dots,x_d^t)$ is independent for infinitely many $t$,
then $\rij xd$  is regular. We are interested in the situation that a (non-\CM)
Noetherian local ring has an independent system of parameters. For instance, the
local ring $\pow K{X,Y}/(X^2,XY)\pow K{X,Y}$ with $K$ a field, does not admit  
an
independent system of parameters, whereas $Y$ is independent in the local ring
$\pow K{X,Y}/(X^2,XY^2)\pow K{X,Y}$.

Before we state this weaker form, we need to introduce one last concept.
Namely, in order to apply Theorem~\ref{T:par}, we have to enforce for the
separated ultraproduct to be reduced and equidimensional. It does not suffice
do require that the components have the same properties. For instance, the
separated
ultraproduct of the analytically irreducible one-dimensional domains $\seq Rw$ 
associated to the cusps $X^2-Y^w$ in $\mathbb A_K^2$ is the the non-reduced curve
$X^2=0$ in $\mathbb A_{\ul K}^2$, where $\ul K$ is the ultraprower of the
field $K$. To control the separated ultraproduct, we use a result of Swanson.
Let us say that a local ring $(A,\pr)$ has \emph{$k$-bounded
multiplication} if for all $n$ and all $a,b\in A$ we have 
$$
\op{ord}(ab) \leq k\max\{\op{ord}(a),\op{ord}(b)\}
$$
where $\op{ord}(a)$ is the supremum of all integers $n$ such that $a\in\pr^n$
(with the usual convention that $\pr^0=A$ and $\op{ord}(a)=\infty$ when $a\in
\inf A$). It is shown in \cite[Theorem 3.4]{Swa} and \cite[Proposition
2.2]{HueSwa} that for a Noetherian
local ring $A$ the multiplication is $k$-bounded for some $k$ \iff\ $A$ is
analytically irreducible (see \cite[Proposition 5.6]{OlbSaySha} for a
characterization in terms of ultraproducts). At present we do not have a good
understanding of the smallest such $k$: in \cite{HueSwa} an upper bound in terms
of Rees valuations is given. It is not hard to see
that having $k$-bounded multiplication is preserved under separated
ultraproducts (see also the next proof). Therefore, the above example on the 
cusps $X^2-Y^w$ shows that an upper bound on $k$ cannot be given in terms of
dimension and multiplicity alone.

\begin{proposition}
For each quadruple $(p,d,l,k)$ with $p$ a prime number, there exists a bound
$N:=N(p,d,l,k)$ with the
following property. Let $R$ be a  $d$-dimensional mixed \ch\  analytically
irreducible local domain  of residual \ch\ $p$. Assume the  multiplication in 
$R$ is $k$-bounded. Let
$\mathbf x$ be a system of parameters generating an ideal of co-length at most $l$. If $\mathbf
xR$ is solidly closed (so that $R$ is weakly S-rational) and the ramification
index of $R$ is at least $N$, then $\mathbf x$ is independent.
\end{proposition}
\begin{proof}
If not, then we get for a fixed quadruple $(p,d,l,k)$ and for each $w\in\nat$ a
counterexample consisting of a $d$-dimensional   mixed \ch\  analytically
irreducible local domain $\seq Rw$ of residual \ch\ $p$, with $k$-bounded
multiplication, and a `dependent' system of parameters $\seq{\mathbf
x}w=(\seq{x_1}w,\dots,\seq {x_d}w)$
in $\seq Rw$ generating
a solidly closed ideal  of co-length at most $l$. This means, after renumbering,
  that there exists $\seq aw\notin\seq{\mathbf x}w\seq Rw$ such that $\seq
aw\seq{x_d}w\in (\seq{x_1}w,\dots,\seq {x_{d-1,}}w)\seq Rw$. Let $a$, $x_i$ and
$R$ be the ultraproducts of the $\seq aw$, $\seq{x_i}w$ and $\seq Rw$
respectively.

By \los, $I:=\rij xdR$ has co-length at most $l$,
 does not contain $a$ and
 is ultra-solidly closed. We leave it
to the reader to verify that $\rij xd$ is a system of cata-parameters (use for
instance the argument in the proof of Theorem~\ref{T:ureg} or \cite[Theorem
8.1]{SchFinEmb}). The multiplication in $R$ is again $k$-bounded by \los,
and it is not hard to show that this implies the same for its separated quotient $\sep
R$. In particular, $\sep R$ is a domain so that we can apply Theorem~\ref{T:par} to conlude that the ideal $I$ is
  cata-tightly closed. On the other hand, \los\ yields that $ax_d\in \rij x{d-1}R$ so that
$a\in \ctc{\rij x{d-1}R}$ by Proposition~\ref{P:CC}, whence $a\in \ctc I=I$,
contradiction.
\end{proof}

Hence combining this result with Theorem~\ref{T:wSrat}, we may add to the
conclusion in the latter theorem that $R$ as above admits an independent system
of parameters.

Our last application gives an `asymptotic' affirmative answer to a question
posed by Hochster concerning solid closure in mixed \ch\ in \cite[Remark
10.13]{HoSol}:

\begin{theorem}
For each quintuple $(p,d,l,m,k)$ with $p$ a prime number, there exists a bound
$N:=N(p,d,l,m,k)$ with the
following property. Let $R$ be a $d$-dimensional mixed \ch\  analytically
irreducible local domain  of residual \ch\ $p$. Assume the  multiplication in 
$R$ is $k$-bounded. Let $\rij
xd$ be a system of parameters generating an ideal of co-length at most $l$ and
let $0\leq a_i<b_i\leq m$.
 If the ramification index of
$R$ is at least $N$, then 
$$
(\sc{(x_1^{b_1},\dots,x_d^{b_d})R} : x_1^{a_1}\cdots x_d^{a_d})\subseteq
\sc{(x_1^{b_1-a_1},\dots,x_d^{b_d-a_d})R}.
$$
\end{theorem}
\begin{proof}
As before, this follows by the same argument from the corresponding
ultra-version, which holds in view of Remark~\ref{R:CC} and
Theorem~\ref{T:par}.   We leave the details to the reader.
\end{proof}

\providecommand{\bysame}{\leavevmode\hbox to3em{\hrulefill}\thinspace}
\providecommand{\MR}{\relax\ifhmode\unskip\space\fi MR }
\providecommand{\MRhref}[2]{%
  \href{http://www.ams.org/mathscinet-getitem?mr=#1}{#2}
}
\providecommand{\href}[2]{#2}


\end{document}